\documentclass[11pt,letterpaper]{article}
\usepackage{fullpage}
\usepackage[margin=1in]{geometry}
\usepackage{multirow}
\usepackage{url}
\usepackage[utf8]{inputenc}
\usepackage{amsmath,amssymb,amsfonts}
\usepackage{graphicx}
\usepackage{color}
\usepackage{fancyhdr}
\usepackage{tikz}
\usepackage{soul}
\usepackage{comment}
\usepackage{bm}
\usepackage{enumitem,kantlipsum}
\usepackage{wrapfig}
\usepackage[compact]{titlesec}
\titlespacing*{\section} {0pt}{1ex}{1ex}
\titlespacing*{\subsection} {0pt}{1ex}{1ex}
\titlespacing*{\subsubsection}{0pt}{1ex}{1ex}

\usepackage{amsopn}

\usepackage{amsthm}
\usepackage[square,sort,comma,numbers]{natbib}
\DeclareMathOperator*{\argmin}{\arg\!\min}

\newcommand{\real}{{\rm{I\hspace{-.75mm}R}}}

\renewcommand*{~}{\relax\ifmmode\sim\else\nobreakspace{}\fi}

\newcommand{\Var}{\mbox{Var}}

\newcommand{\Fbar}{\bar{F}}

\newcommand{\Fhat}{\hat{F}}

\newcommand{\rhohat}{\hat{\rho}}
\newcommand{\sigmahat}{\hat{\sigma}}

\newcommand{\Ftilde}{\widetilde{F}}

\newcommand{\BFc}{\bm{c}}

\newcommand{\BFx}{\bm{x}}

\newcommand{\BFX}{\bm{X}}

\newcommand{\BFDelta}{\bm{\Delta}}

\newcommand{\mcB}{\mathcal{B}}

\newcommand{\mcO}{\mathcal{O}}

\newcommand{\mbE}{\mathbb{E}}

\newcommand{\mbN}{\mathbb{N}}

\usepackage{array, caption, floatrow, tabularx, makecell, booktabs}%
\setcellgapes{3pt}
\usepackage{algorithmic}
\usepackage{algorithm}

\newtheorem{theorem}{Theorem}

\newtheorem{lemma}{Lemma}

\newtheorem{proposition}{Proposition}

\newtheorem{corollary}[theorem]{Corollary}

\newtheorem{definition}{Definition}

\newtheorem{assumption}{Assumption}

\newtheorem{remark}{Remark}

\usepackage{bm}
\usepackage{multicol}
\usepackage{lipsum}
\usepackage{mwe}
\usepackage{graphicx}
\usepackage{subfig}
\usepackage{amssymb}
\usepackage{booktabs}
\usepackage{enumitem}
\usepackage{cleveref}

\usepackage{multirow}
\usepackage{amsfonts}
\usepackage{graphicx}
\usepackage{epstopdf}
\usepackage{algorithmic}
\usepackage{algorithm}

\usepackage{qcircuit}
\usepackage{physics}

\usepackage{authblk}

\ifpdf
  \DeclareGraphicsExtensions{.eps,.pdf,.png,.jpg}
\else
  \DeclareGraphicsExtensions{.eps}
\fi

\makeatletter
\newcommand{\ALOOP}[1]{\ALC@it\algorithmicloop\ #1%
  \begin{ALC@loop}}
\newcommand{\ENDALOOP}{\end{ALC@loop}\ALC@it\algorithmicendloop}

\newcommand{\algorithmicbreak}{\textbf{break}}
\newcommand{\BREAK}{\STATE \algorithmicbreak}
\makeatother

\title{Multi-Fidelity Stochastic Trust Region Method with Adaptive Sampling}
\author[1]{Yunsoo Ha\thanks{yunsoo.ha@nrel.gov}}
\author[1]{Juliane Mueller\thanks{juliane.mueller@nrel.gov}}
\affil[1]{Computational Science Center, National Renewable Energy Laboratory, \protect\\ 15013 Denver West Parkway, Golden, 80401, Colorado, USA}
\date{}

\begin{document}
\vspace{-3 mm}

\maketitle

\section*{ABSTRACT}
Simulation optimization is often hindered by the high cost of running simulations. Multi-fidelity methods offer a promising solution by incorporating cheaper, lower-fidelity simulations to reduce computational time. However, the bias in low-fidelity models can mislead the search, potentially steering solutions away from the high-fidelity optimum. To overcome this, we propose ASTRO-MFDF, an adaptive sampling trust-region method for multi-fidelity simulation optimization. ASTRO-MFDF features two key strategies: (i) it adaptively determines the sample size and selects { appropriate sampling strategies to reduce computational cost}; and (ii) it selectively uses low-fidelity information only when a high correlation with the high-fidelity is anticipated, reducing the risk of bias. We validate the performance and computational efficiency of ASTRO-MFDF through numerical experiments using the SimOpt library.

\section{INTRODUCTION}
\label{sec:intro}
Simulation optimization (SO) has become a key method for optimizing objective functions in uncertain environments, gaining significant attention for its ability to tackle real-world problems involving randomness and complex systems, such as quantum computing and renewable energy \cite{ha2025two,sakki2022renewable}. However, its practical implementation can be challenging due to high computational cost of evaluating a stochastic function value. To address this, simulation { models } can be { developed at different levels of fidelity, which sometimes requires substantial modeling or coding effort, } a method known as multi-fidelity (MF) simulation. In MF simulation, models are designed with a hierarchical structure, where high-fidelity models provide detailed and accurate representations of the process, while low-fidelity models offer a more computationally efficient, simplified version. For example, a high-fidelity model may represent the full manufacturing process in detail, whereas lower-fidelity models might omit certain machines that are not critical \cite{zhang2022improved}.
Another way to construct MF simulation models is by modifying the length of the simulation runs. When aiming to optimize systems under a steady-state condition, using shorter run lengths results in lower-fidelity models that produce less accurate output estimates but require reduced computational effort \cite{chen2017stochastic}.

In this paper, we consider the multi-fidelity simulation optimization (MFSO) problem 
\begin{equation}
    \min_{\BFx\in \real^d} f^0(\BFx) := \mbE_{\Xi^0}[F^0(\BFx,\xi^0)],
    \label{eq:problem}
\end{equation}
where $f^0:\real^d \rightarrow \real$ is nonconvex and has a lower bound, $F^0:\real^d\times\Xi^0 \rightarrow \real$ is a random function, and $\xi^0:\Omega\rightarrow\Xi^0$ is a random element. Here, the index $0$ represents the highest-fidelity simulation, while increasing index values correspond to lower-fidelity simulations. We consider zeroth-order stochastic oracles, where the derivative information is not directly available from the Monte Carlo simulation. We allow each realization $F^0(\BFx,\xi^0_i)$ to be nonconvex and nonsmooth, as is common in complex simulation models.
Since we only have access to $F^0(\BFx)$, we can estimate $f^0(\BFx)$ by $\Fbar^0(\BFx,n)=n^{-1}\sum_{i=1}^n F^0(\BFx,\xi^0_i)$, variance of $F^0(\BFx,\xi^0)$ by $(\sigmahat^0(\BFx,n))^2=n^{-1}\sum_{i=1}^n(F^0(\BFx,\xi^0_i)-\Fbar^0(\BFx,n))^2$, and covariance between $F^0(\BFx,\xi^0)$ and $F^i(\BFx,\xi^i)$ by $\sigmahat^{0,i}(\BFx,n) = n^{-1}\sum_{j=1}^n(F^0(\BFx,\xi^0_j)-\Fbar^0(\BFx,n))(F^i(\BFx,\xi^i_j)-\Fbar^i(\BFx,n))$ for any $i \in \mbN.$ Typically, iterative algorithms are used to solve \eqref{eq:problem} by generating a random sequence of iterates  $\{\BFX_k\}$ that converges toward an optimal solution. In the context of MFSO, the goal is to accelerate this convergence by effectively leveraging information from MF simulations.


In deterministic MF problems, multi-fidelity Bayesian optimization (MFBO) is widely used, with co-kriging serving as the surrogate model \cite{do2023multi}. However, applying MFBO to solve \eqref{eq:problem} requires two key assumptions \cite{foumani2023effects}. First, function estimates must be sufficiently accurate. Second, lower-fidelity functions should exhibit a strong correlation with the high-fidelity function throughout the entire search space. Otherwise, lower-fidelity data can negatively impact the co-kriging model, slowing down optimization even more than using only the high-fidelity model. To mitigate these challenges, it is crucial to carefully select both the sample size $(n)$ and the regions where lower-fidelity data provides meaningful information. This consideration leads to the following research question:
\begin{enumerate}
\centering
\item[(RQ)] \textit{Where, which fidelity, and how many times should oracles be queried to efficiently solve \eqref{eq:problem}?}
\end{enumerate}

The adaptive sampling trust region method (ASTRO-DF) is one of the most effective algorithms for addressing this question, as it dynamically adjusts the sample size and trust region \cite{Sara2018ASTRO}. Specifically, the sample size is chosen to balance the trade-off between the optimality gap and estimation errors, helping to reduce the computational burden. In addition, the search space is confined to a neighborhood of the current iterate, known as the trust region, by constructing and optimizing a surrogate model within this region. The trust region is updated based on whether the surrogate model provides a sufficiently accurate approximation of the true objective function, enabling correlation measurement of the MF functions within a localized search space. See \cite{katya:DFObook} for the details of deterministic trust region method for derivative-free optimization. To address (RQ) in the bi-fidelity setting, we recently proposed an extension of ASTRO-DF, known as ASTRO-BFDF \cite{ha2024adaptive}. ASTRO-BFDF differs from ASTRO-DF in two primary aspects: i) it maintains two separate trust regions for high-fidelity and low-fidelity functions, and ii) it employs an adaptive sampling approach that utilizes both Bi-fidelity Monte Carlo (BFMC) and { standard } Monte Carlo (MC) to reduce the computational burden. 

BFMC is a special case of the multi-fidelity Monte Carlo (MFMC) method that uses exactly two levels of simulation fidelity. MFMC has been introduced in~\cite{karen2016mfmc} to enhance the efficiency of Monte Carlo methods by leveraging MF simulation oracles through a control variate approach. Specifically, the MFMC estimate is obtained by 
\begin{equation} \label{eq:MFMC}
    \Fhat^t(\BFx,\bm{n},\bm{c}) = \frac{1}{n^t}\sum_{j=1}^{n^t} F^t(\BFx,\xi^t_j) + \sum_{i=t+1}^{q} c^i \left(\frac{1}{n^i}\sum_{j=1}^{n^{i}} F^i(\BFx,\xi^i_j) - \frac{1}{n^{i-1}}\sum_{j=1}^{n^{i-1}} F^i(\BFx,\xi^i_j) \right),
\end{equation}
where $q$ denotes the index of the lowest-fidelity simulation, $t$ is the target fidelity level, $\bm{c}=\{c^1,c^2,\cdots,c^q\}$ with $c^i$ indicating the coefficient associated with $i$-th fidelity oracle, and $\bm{n} =\{n^0, n^1, \dots, n^q\}$ with $n^i$ denoting the sample size used for the $i$-th fidelity oracle. Given fixed $\bm{n}$ and $\bm{c}$ such that $n^0 < n^1 < \dots < n^q$, the MFMC estimate is unbiased estimate for $f^0(\BFx)$ with a variance
\begin{equation} \label{eq:variance-mfmc-main}
        \Var(\Fhat^t(\BFx,\bm{n},\bm{c})) = \frac{(\sigma^t(\BFx))^2}{n} + \sum_{i=t+1}^q\left(  \frac{1}{n^{i-1}}-\frac{1}{n^i}\right) \left( (c^i)^2(\sigma^i(\BFx))^2 - 2c^i \sigma^{0,i}(\BFx)\right),
\end{equation}
where $(\sigma^{i}(\BFx))^2$ is the variance of $F^i(\BFx,\xi^i)$ and $\sigma^{0,i}$ is the covariance between $F^0(\BFx,\xi^0)$ and $F^i(\BFx,\xi^i)$. Suppose that we set $t=0$ and $q=1$. Then the variance reduction { is } achieved when the condition $2c^1 \sigma^{0,1}(\BFx) \ge (c^1)^2(\sigma^1(\BFx))^2$ holds. {To satisfy this condition, we employ the Common Random Numbers (CRN) technique, which facilitates positive correlation between fidelity levels, and later propose a principled choice of $\bm{c}$ to minimize the variance of the MFMC estimator.}


In this paper, we extend ASTRO-BFDF to address the MFSO problem \eqref{eq:problem} by incorporating the MFMC technique, which we refer to as ASTRO-MFDF. We first present a multi-fidelity adaptive sampling algorithm in Section \ref{sec:MFAS}. Then, we introduce the ASTRO-MFDF algorithm in Section \ref{sec:ASTRO-MFDF}. Finally, we demonstrate the performance of ASTRO-MFDF through experiments with stochastic Rosenbrock functions and continuous $(s,S)$ inventory problems in Section \ref{sec:numerical}. Throughout the paper, we use capital letters for random variables and bold font for vectors.

\section{Multi-Fidelity Adaptive Sampling}
\label{sec:MFAS}
The adaptive sampling method used within ASTRO-DF performs sequential estimation to effectively reduce computational effort associated with MC. Specifically, the sample size at $\BFx \in \real^d$ has been determined by
\begin{equation}
\label{eq:ori-AS}
    N_k(\BFx) = \min\biggl\{ n \in \mbN: \underbrace{\frac{\max\{\sigma_{lb},\sigmahat^0(\BFx,n)\}}{\sqrt{n}}}_{{\text{the stochastic error}}}\leq \underbrace{\frac{\kappa\Delta_k^{2}}{\sqrt{\lambda_k}}}_{{\text{the optimality gap}}}\biggr\},
\end{equation}
where $\kappa$, $\lambda_k$, and $\sigma_{lb}$ are positive constants, $\Delta_k$ is the trust region size at iteration $k$ \cite{ha2023}. The condition in \eqref{eq:ori-AS} balances the stochastic error and the optimality gap, represented by the standard deviation estimate and the square of the trust region size respectively. Following the same principle, when employing the MFMC estimate, our goal is to determine $\bm{n}$ and $\bm{c}$ such that the condition in \eqref{eq:ori-AS} is satisfied by replacing {the stochastic error term under MC, i.e., } $n^{-1}(\sigmahat^0(\BFx,\bm{n}))^2$, with { that under MFMC, i.e., } the estimate of $\text{Var}(\hat{F}^{0}(\BFx,\bm{n},\bm{c}))$. However, incorporating the adaptive sampling strategy to MFMC introduces the following two challenges:
\begin{enumerate}
    \item[(C1)] The decision variables become multi-dimensional, consisting of the vectors $\bm{n}$ and $\bm{c}$. In the MC case, the decision variable is one-dimensional, specifically the scalar $n$. Thus, to determine $N_k(\BFx)$, we initially compute $\sigmahat^0(\BFx,n=3)$ and increment $n$ by one until the condition in \eqref{eq:ori-AS} is satisfied, due to the unknown variance $(\sigma^0(\BFx))^2$. Similarly, the variances and covariances appearing in \eqref{eq:variance-mfmc-main} are also unknown, necessitating sequential estimation for these quantities as well. However, designing a straightforward sequential estimation algorithm to determine $\bm{n}$ and $\bm{c}$ while minimizing computational burden is challenging.
    \item[(C2)] Considering factors such as the variances and covariances presented in \eqref{eq:variance-mfmc-main} and the querying costs of the MF simulation oracles, the computational cost of MFMC can be higher than that of MC for achieving the same level of accuracy, despite MFMC being a variance reduction method. { Hence, } we must determine which MC method (MC or MFMC) to use while sequentially estimating variances and covariances, introducing an additional binary decision variable. 
\end{enumerate}

In (C1), we should determine the optimal $\bm{n}$ and $\bm{c}$ based on the available information. Suppose, for example, that variance and covariance estimates are currently available with sample sizes $\widetilde{\bm{n}}^i$ for the $i$-th fidelity simulation for any $\widetilde{\bm{n}}^i \ge 2$. Using these estimates, we can obtain an approximation of \eqref{eq:variance-mfmc-main}, specifically, $\widehat{\Var}(\hat{F}^{0}(\BFx,\bm{n},\BFc))$, by replacing $\sigma^i(\BFx)$ and $\sigma^{0,i}(\BFx)$ with their corresponding estimators $\widehat{\sigma}^i(\BFx,\widetilde{n}^i)$ and $\widehat{\sigma}^{0,i}(\BFx,\widetilde{n}^0)$, respectively. We can then determine the optimal solutions for $\bm{n}$ and $\bm{c}$, denoted by $\bm{n}_*$ and $\bm{c}_*$, by solving the following optimization problem:
\begin{equation} \label{eq:mfmc-problem}
\begin{split}
    [\bm{n}_*,\bm{c}_*]  \in \argmin_{\bm{n},\bm{c} \in \real^q} \quad \sum_{i=0}^q w^i n^i  \\
    \text{subject to} \quad
     \widehat{\Var}(\hat{F}^{0}(\BFx,\bm{n},\bm{c})) &\le \kappa^2 \Delta_k^4 \lambda_k^{-1} \\
     n^i - n^{i+1} &\le 0 \quad \forall i\in \{0,1,\dots,q-1\}\\
     \widetilde{n}^i - n^i  &\le 0 \quad \forall i\in \{0,1,\dots,q\},
\end{split}
\end{equation}
where $\bm{w} = \{w^0,w^1,\dots,w^q\}$ with $w^i$ indicating the cost of querying the $i$-th fidelity simulation oracle. The first constraint originates from the condition stated in \eqref{eq:ori-AS}. The predicted optimal cost of the MFMC method then becomes $\sum_{i=0}^q w^i n_*^i$. Additionally, the predicted optimal cost for the MC method is calculated as $w^0 n_p$, where $n_p = \lceil (\sigmahat^0(\BFx,\widetilde{n}^0))^2 \lambda_k (\kappa^2 \Delta_k^4)^{-1} \rceil$, based on the condition from \eqref{eq:ori-AS}. Therefore, MFMC is selected if $\sum_{i=0}^q w^i n_*^i \le w^0 n_p^0$; otherwise, MC is chosen. This resolves the issue mentioned in (C2). {It is worth noting that problem \eqref{eq:mfmc-problem} is a nonlinear continuous optimization problem with mostly linear constraints and a single nonlinear constraint, where the problem dimension is relatively small. Since the primary computational cost in simulation optimization lies in the simulation itself, the cost of solving problem \eqref{eq:mfmc-problem} is assumed to be negligible. We also note that the optimal integer sample sizes can be obtained by rounding up $\bm{n}_*$.}

\begin{figure} [htp]
\centering
\includegraphics[width=1\columnwidth]{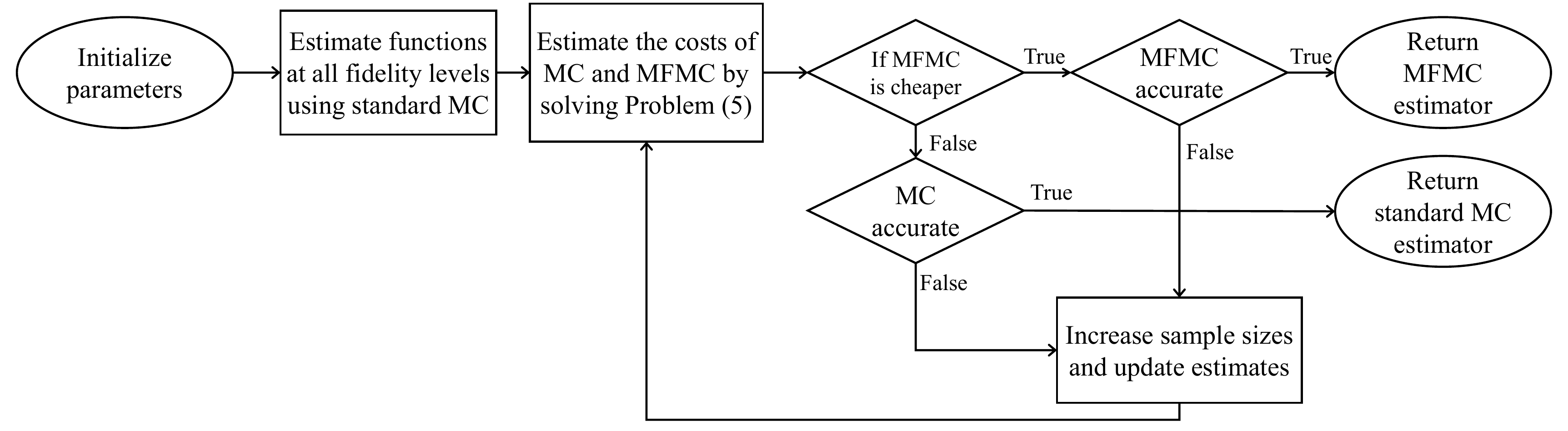}
\caption{Flowchart of the MFAS procedure. The algorithm first initializes relevant parameters and estimates function values at all fidelity levels using standard MC sampling. It then compares the estimated costs and variances of standard MC and MFMC estimators by solving Problem \eqref{eq:mfmc-problem}. Depending on the cost and accuracy, the algorithm either returns a selected estimator or increases the sample sizes and updates the estimates adaptively.}
\label{fig:mfas}
\end{figure}

If the chosen estimate is sufficiently accurate—specifically, if it satisfies the condition in \eqref{eq:ori-AS} based on its own variance estimate—then that function estimate can directly serve as the final output of the multi-fidelity adaptive sampling method. If this accuracy condition is not met, further replications are required to refine the estimate. In this case, when MC has been chosen, the $0$-th fidelity oracle naturally becomes the default choice. 
Conversely, if MFMC is more economical but the accuracy condition $\widehat{\Var}(\hat{F}^{0}(\BFx,\widetilde{\bm{n}},\bm{c}_*)) > \kappa^2 \Delta_k^4 \lambda_k^{-1}$ remains unmet, {i.e., $\widetilde{n}^j < n_*^{j}-1$ for some $j\in\{0,1,\dots,q\}$, we proceed by querying the highest such fidelity level, as allocating additional samples to higher fidelity levels is generally more effective in reducing the variance. } 
After selecting the appropriate fidelity level, we conduct additional replications, update variance and covariance estimates accordingly, and then resolve the optimization problem \eqref{eq:mfmc-problem}. This iterative process continues until one of the accuracy conditions is ultimately satisfied.
Although we estimated $f^0$ above, target fidelity can be any $i \in \{0,1,\dots,q\}$. To complement the algorithmic description, we include a flowchart that illustrates the overall logic of the multi-fidelity adaptive sampling method (MFAS), including the decision-making process between MC and MFMC, as well as the adaptive sample size update based on accuracy checks (see Figure~\ref{fig:mfas}).

\section{Multi-fidelity stochastic Trust Region Method with MFAS}
\label{sec:ASTRO-MFDF}

\begin{figure} [htp]
\centering
\includegraphics[width=1\columnwidth]{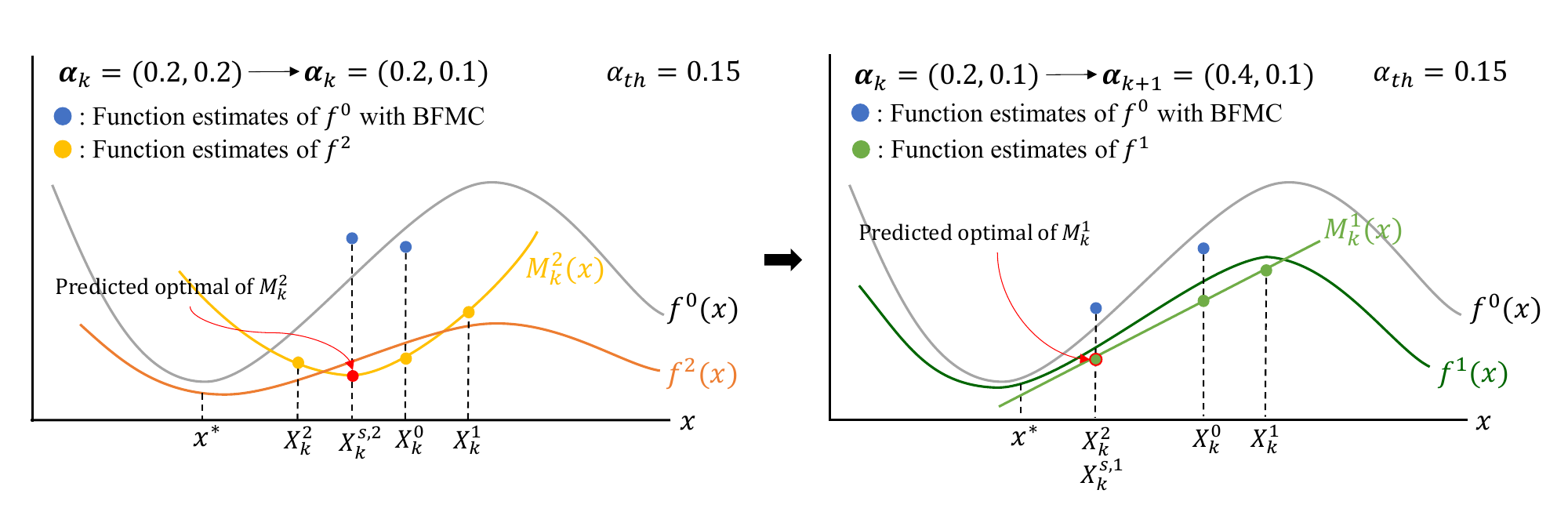}
\caption{Illustration of the inner loop (Steps 2–9) of ASTRO-MFDF, which iteratively applies Algorithm \ref{alg:TRO-LFDF} across fidelity levels. Since $M_k^2$ fails to generate a better candidate than $\BFX_k^0$, $\alpha_k^2$ decreases, and $M_k^1$ is instead utilized, successfully yielding a better candidate. As a result, at iteration $k+1$, $M_k^1$ will be constructed first based on the updated $\bm{\alpha}_{k+1}$.}
\label{fig:astromf}
\end{figure}

\begin{figure} [htp]
\centering
\includegraphics[width=0.5\columnwidth]{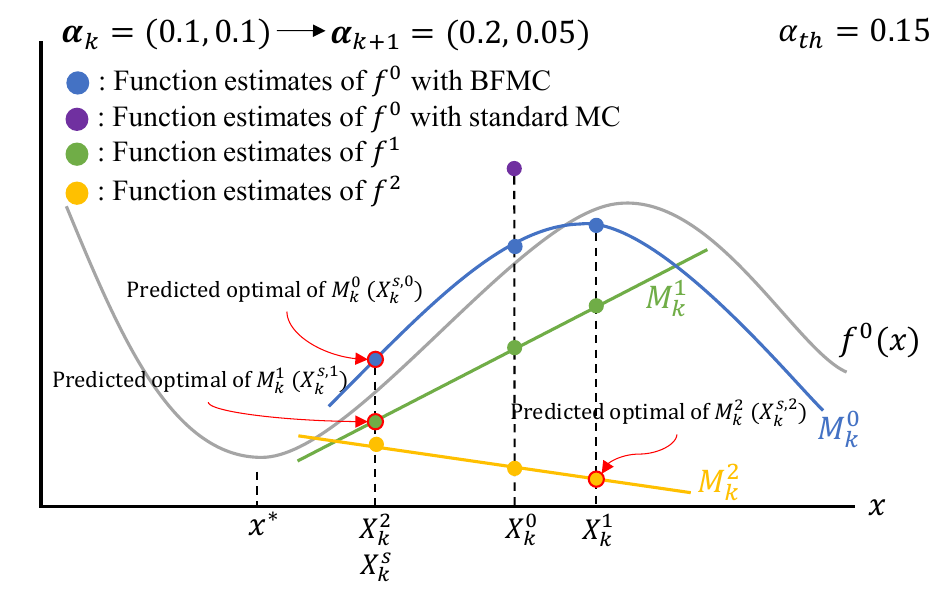}
\caption{{Illustration of Steps 10–15 in ASTRO-MFDF. Since all $\alpha_k^i < \alpha_{\text{th}}$, design points are shared across all models. Lower-fidelity estimates (green and yellow dots) are used in the BFMC estimator (blue dots) to approximate $f^0$ and are then reused in constructing $M_k^1$ and $M_k^2$. As $M_k^2$ fails to improve over the incumbent, $\alpha_k^2$ decreases, while $M_k^1$ succeeds, increasing $\alpha_k^1$.}}
\label{fig:fail}
\end{figure}

\begin{algorithm}[htp]  
\caption{\texttt{ASTRO-MFDF}}
\label{alg:TRO-MFDF}
\begin{algorithmic}[1]
\REQUIRE Initial incumbent $\BFx_{0}\in\real^d$, initial and maximum trust region radius $\BFDelta_0$ and $\Delta_{\max}>0$, model fitness thresholds $0<\eta<1$ and certification threshold $\mu>0$, expansion and shrinkage constants $\gamma_1>1$ and $\gamma_2\in(0,1)$, sample size lower bound sequence $\{\lambda_k\} = \{\mcO(\log k)\}$, adaptive sampling constant $\kappa>0$, correlation vector $\bm{\alpha}_0=\{\alpha^1_0,\alpha^2_0,\dots,\alpha^t_0\}$, and lower bound of an initial variance approximation $\sigma_{lb}>0$, and sufficient reduction constant $\zeta>0$. 
\FOR{$k=0,1,2,\hdots$}
\FOR{$t=q,q-1,\dots,1$}
\STATE \label{ASBFTRO:callingLF} Obtain $\BFX_k^s$ and $\Delta_k^t$ by calling Algorithm \ref{alg:TRO-LFDF} (\texttt{ASTRO-LFDF-}$t$).
\IF{{Algorithm \ref{alg:TRO-LFDF} yielded a candidate with sufficient HF function reduction (i.e., $\BFX_k^s \not=\BFX_k^0$)}}
\STATE Set $(\BFX_{k+1},\Delta^{t}_{k+1})= (\BFX_k^s, \gamma_1\Delta^{t}_k)$ and $\alpha^t_{{k+1}} = \gamma_1 \alpha^t_k$.
\STATE Set $\Delta_{k+1}^j = \max\{\Delta_k^{j}, \Delta_{k+1}^t\}$ for all $j \in \{0,1,\dots,t-1\}$ and $k = k+1$.
\BREAK
\ENDIF
\ENDFOR
\IF{{Algorithm \ref{alg:TRO-LFDF} failed to find a better candidate than the current solution (i.e., $\BFX_k^s =\BFX_k^0$)}}
\STATE \label{ASBFTRO:designsetselect} Select $\{ \BFX_{k}^{i}\}_{i=0}^{2d}\subset\mcB(\BFX_{k};\Delta^0_{k})$.

\STATE \label{ASBFTRO:estimate-hf} Estimate the $t$-fidelity function $\Ftilde^t$ at $\{ \BFX_{k}^{i}\}_{i=0}^{2d}$ using BFAS with $\Delta_k = \Delta_k^0$, construct local models $M_k^{t}(\BFX)$, approximately compute the minimizers $\BFX^{s,t}_k\in\argmin_{\left\| \BFX-\BFX_k \right\| \leq \Delta^t_{k}}M^t_{k}(\BFX)$ and then, estimate the $0$-fidelity function $\Ftilde^0(\BFX^{s,t}_k)$ using BFAS with $\Delta_k = \Delta_k^0$ for all $t\in \{0,1,\dots,q\}$.

\STATE Set the candidate point $\BFX^{s}_k \in \argmin_{\BFx\in\{\BFX^{s,i}_k\}^q_{i=0}} \Ftilde^0(\BFx),$ and compute the success ratio $\rhohat_k$ and $\rhohat^{t}_k$ for any $t \in \{1,2,\dots,q\}$ as
\begin{equation*}
    \rhohat_k = \frac{\Ftilde^0(\BFX_k^0)-\Ftilde^0(\BFX_k^s)}{M_k^0(\BFX_k^0)-M_k^0(\BFX_k^{s,{0}})}
    \text{ and }    
    \rhohat^{t}_k = \frac{\Ftilde^0(\BFX_k^0)-\Ftilde^0(\BFX_k^{s,t})}{\max\{\zeta(\Delta_k^0)^2,M_k^{t}(\BFX_k^0)-M_k^{t}(\BFX_k^{s,t})\}}.
\end{equation*}
\STATE If $M_k^t$ succeeds (i.e., $\rhohat^{t}_k \ge \eta)$, set ${\alpha^t_{k+1}}=\gamma_1 \alpha^t_k$; otherwise set ${\alpha^t_{k+1}}=\gamma_2 \alpha^t_k$ for all $t \in \{1,2,\dots,q\}$.
\STATE \label{HF:delta-l-update} 
Update $(\BFX_{k+1},\Delta^0_{k+1})$ as follows:
\[
(\BFX_{k+1},\Delta^0_{k+1}) =
\begin{cases}
    (\BFX_k^s,\min\{\gamma_1\Delta^0_{k}, \Delta_{\max}\}) & \text{if } \rhohat_k \ge \eta \text{ and } \mu\|\nabla M^0_{k}(\BFX_{k})\| \ge \Delta^0_{k}, \\
    (\BFX_{k},\gamma_2\Delta^0_{k}) & \text{otherwise}.
\end{cases}
\]
\STATE Set $\Delta^{t}_{{k+1}} = \min\left\{\Delta^{t}_k,\Delta^0_k\right\}$ for all $t \in \{1,2,\dots,q\}$, and $k=k+1$.
\ENDIF
\ENDFOR
\end{algorithmic}
\end{algorithm}

\begin{figure} [htp]
\centering
\includegraphics[width=1\columnwidth]{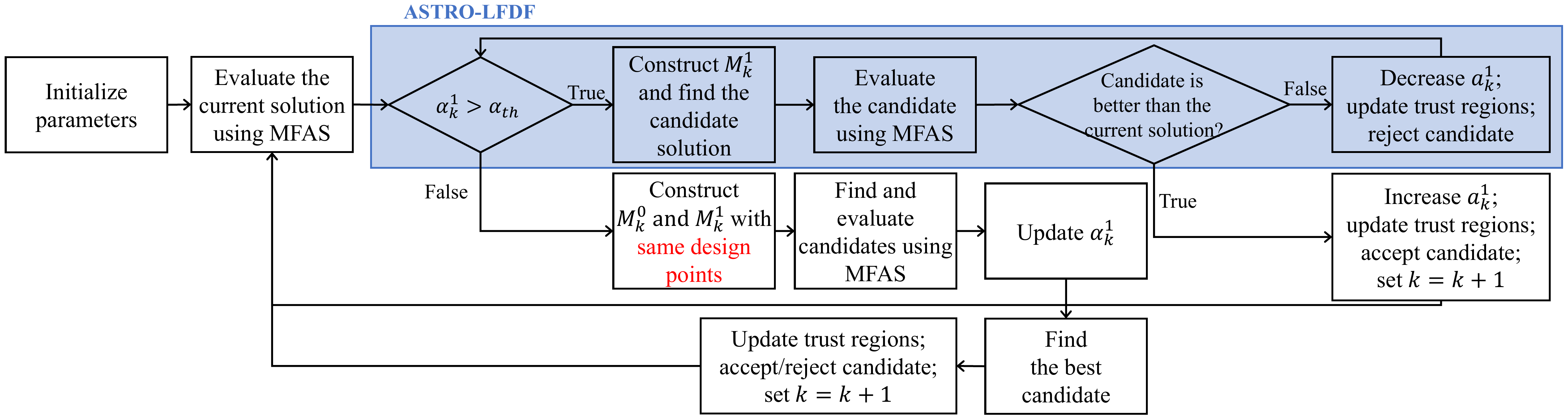}
\caption{Flowchart of ASTRO-MFDF, illustrating the integration of ASTRO-LFDF (highlighted in blue), the decision process based on the correlation measure, the generation of candidate solutions, and the update of trust regions and correlation values. For clarity, the flowchart depicts the bi-fidelity case ($q=1$) as a special instance.}
\label{fig:astro-mfdf-flowchart}
\end{figure}

\noindent
In traditional stochastic trust region methods \cite{STRONG,chen2018storm}, a single interpolation/regression model ($M_k$) is constructed at iteration $k$, and the next candidate point is determined by { approximately } minimizing this model within the trust region, {usually by computing a Cauchy point or using iterative methods such as the conjugate gradient method; see \cite{Nocedal:NumericalOptbook} for further details. }
The key distinction of ASTRO-MFDF is that, instead of relying on a single model, it can construct multiple interpolation models, denoted as $M_k^i$ for any $i \in \{0,1,\dots,q\}$, across multiple trust regions, represented by 
$\BFDelta_k:=\{\Delta_k^0,\Delta_k^1, \dots, \Delta_k^q\}$. The local model constructed with lower-fidelity oracles is prioritized in selecting the next iterate, effectively guiding the iterates closer to the optimal solution using only lower-fidelity simulations. However, constructing $q$ local models at every iteration can be cumbersome and a waste of computational resources, especially since lower-fidelity simulations may not yield better solutions in certain feasible regions due to the inherent bias between the high-fidelity function and the lower-fidelity functions. Hence, we employ an adaptive correlation vector $\bm{\alpha}_k := \{\alpha_k^1,\dots,\alpha_k^q\}$, updated dynamically to capture the correlation between $f^0$ and $f^i$ for any $i\in\{1,2,\dots,q\}$. Specifically, if the candidate point generated by $M_k^i$ achieves a sufficient reduction in the estimated objective function, $\alpha_k^i$ increases; otherwise, it decreases. When $\alpha_k^i$ exceeds a threshold $\alpha_{th}$ (See Figure \ref{fig:astromf}), $M_k^i$ is constructed  using the individual design set, as a strong correlation between $f^0$ and $f^i$ is expected. If $\alpha_k^i < \alpha_{th}$ for all $i \in \{1,2,\dots,q\}$ (See Figure \ref{fig:fail}), the main local model defaults to $M_k^0$ for finding the candidate point. Meanwhile, $M_k^i$ is constructed using same design points for $0$-th fidelity simulation with the main purpose of updating $\alpha_k^i$, even though some design points may lie outside $\mcB(\BFX_k,\Delta_k^i)$, the trust region centered at $\BFX_k$ with radius $\Delta_k^i$. It is important to note that constructing $M_k^i$ with same design points requires minimal additional computational effort, as numerous replications of the $i$-th fidelity simulation are already available from the MFMC estimates in MFAS during the construction of $M_k^0$. To provide an overview of the ASTRO-MFDF procedure, we include a flowchart summarizing the main logic, including the role of the ASTRO-LFDF subroutine, the adaptive update of $\bm{\alpha}_k$, and the candidate selection based on model performance (Figure \ref{fig:astro-mfdf-flowchart}). While the flowchart depicts the bi-fidelity setting as a special case, the logic extends naturally to the general multi-fidelity setting. The pseudo code of ASTRO-MFDF is listed in Algorithm \ref{alg:TRO-MFDF}. 

We now outline some minor details of ASTRO-MFDF for practical purposes. First, we ensure that $\Delta_k^i \ge \Delta_k^j$ for any $i < j$. In practice, the main computational challenge arises from the increasing sample size, which scales at a rate of $\Delta_k^{-4}$ as $\Delta_k$ converges to zero. Therefore, maintaining a large trust region for higher-fidelity functions while progressing toward better iterates with a relatively smaller trust region for lower-fidelity functions enhances finite-time efficiency. Second, when low correlation between $f^0$ and $f^i$ for all $i \in \{1,2,\dots,q\}$ is predicted, i.e., $\alpha_k^i < \alpha_{th}$ for all $i \in \{1,2,\dots,q\}$, we generate $q+1$ candidate points ($\BFX_k^{s,q}$) by minimizing $M_k^i$ within $\Delta_k^0$ for each $i\in\{0,1,\dots,q\}$ (Step \ref{ASBFTRO:estimate-hf} in Algorithm \ref{alg:TRO-MFDF}). The final candidate ($\BFX_k^s$) is then selected as the point with the lowest function estimate, maximizing the likelihood of a successful iteration. Lastly, for an iteration with a candidate point from $M_k^t$ for any $t \in \{1,2,\dots,q\}$ to be considered successful, we impose the additional condition $\Ftilde^0(\BFX_k^0)-\Ftilde^0(\BFX_k^{s,q})\ge \zeta \eta (\Delta_k^0)^2$ to prevent success arising from a negligible model reduction $M_k^t(\BFX_k^0)-M_k^t(\BFX_k^{s,t})$ (Step \ref{LFDF:success-ratio} in Algorithm \ref{alg:TRO-LFDF}).

\begin{algorithm}[htp]  
\caption{[$\BFX_k^s$, $\Delta_k^t$] = \texttt{ASTRO-LFDF-}$t$}
\label{alg:TRO-LFDF}
\begin{algorithmic}[1]
\REQUIRE $\BFX_k^0$, $\Delta^t_{k}$, model fitness thresholds $0<\eta<1$, expansion and shrinkage constants $\gamma_1>1$ and $\gamma_2\in(0,1)$, sample size lower bound sequence $\{\lambda_k\}=\{\mcO(\log k)\}$, adaptive sampling constant $\kappa>0$, correlation constant $\alpha_k^t > 0$, correlation threshold $\alpha_{th}>0$, lower bound of an initial variance approximation $\sigma_{lb}>0$, and sufficient reduction constant $\zeta>0$. 

\ALOOP{}
\IF{$\alpha^t_k < \alpha_{th}$}
\STATE Set $\BFX_k^{s,t} = \BFX_k^0$
\BREAK
\ENDIF
\STATE \label{ASLFTRO:designsetselect} Select the design set $\{ \BFX_{k}^{i}\}_{i=0}^{2d}\subset\mcB(\BFX_{k};\Delta^{t}_{k})$.
\STATE Estimate $t$-fidelity function $\Ftilde^{t}(\BFX_k^i)$ at $\{ \BFX_{k}^{i}\}_{i=0}^{2d}$ using MFAS with $\Delta_k = \Delta_k^t$, construct local model $M_k^{t}(\BFX)$, and approximately compute $\BFX^{s,t}_{k}=\argmin_{\left\| \BFX-\BFX_k \right\| \leq \Delta^{\ell}_{k}}M^{\ell}_{k}(\BFX).$

\STATE Estimate the $0$-fidelity function $\Ftilde^0_k(\BFX_k^{s,t})$ and $\Ftilde^0_k(\BFX_k^0)$ using MFAS with $\Delta_k = \Delta_k^{t}$.
\STATE \label{LFDF:success-ratio} Compute the success ratio $\rhohat_k = \Ftilde_k^0(\BFX_k^0)-\Ftilde_k^0(\BFX_k^{s,t})/\max\{\zeta(\Delta_k^0)^2,M_k^{t}(\BFX_k^0)-M_k^{t}(\BFX_k^{s,t})\}.$
\IF{$\rhohat_k \ge \eta$} \label{LF:sufficient-test} 
\BREAK
\ENDIF
\STATE Set $\Delta_k^{t} = \gamma_2\Delta_k^{t}$ and $\alpha_k^t = \gamma_2 \alpha_k^t$
\ENDALOOP

\RETURN [{$\BFX_k^{s} = \BFX_k^{s,t}$}, $\Delta_k^t$]

\end{algorithmic}
\end{algorithm}

\section{Numerical Experiments}
\label{sec:numerical}
In this section, we analyze the finite-time performance of solvers such as ASTRO-MFDF, ASTRO-DF, and Nelder-Mead, across two different problems. First, we test the solvers on a stochastic Rosenbrock function. Second, we use continuous $(s, S)$ inventory problems to evaluate performance in steady-state simulation optimization, where $s$ is a threshold for ordering and $S$ is a maximum inventory level.

The experiments in this study were carried out using SimOpt, a benchmarking framework designed for evaluating simulation optimization algorithms \cite{eckman2023simopt2}. The experimental process consists of two stages. In the first stage, we generate 20 independent optimization runs (macro-replications) for each solver on each problem, where each run represents a distinct stochastic trajectory. Every solver is allocated a fixed budget specific to the problem at hand and must strategically decide how to utilize this budget, including selecting locations, fidelity levels, and sample sizes. For instance, if the budget is set to a specific number, such as $B$, it allows for $B$ queries to the simulation oracle to solve the given problem. In the MF scenario, we introduce a cost ratio vector $(\bm{w})$, whose component satisfies $w^0=1$ for the highest fidelity and decrease thereafter, i.e., $w^0 > w^1 > \dots$. Therefore, if we query 10 oracles at fidelity level 0 and 20 oracles at fidelity level 1, the total budget consumed is calculated as $10w^0 + 20w^1$. Throughout each run, the solver estimates objective values at various solutions (design points), using a chosen number of replications to guide its search. The second stage involves evaluating the recommended solutions from each run using 200 additional simulations to estimate the true performance of the solutions.

\subsection{Experiments on stochastic Rosenbrock function}
\begin{figure} [htp]
\centering
\includegraphics[width=1\columnwidth]{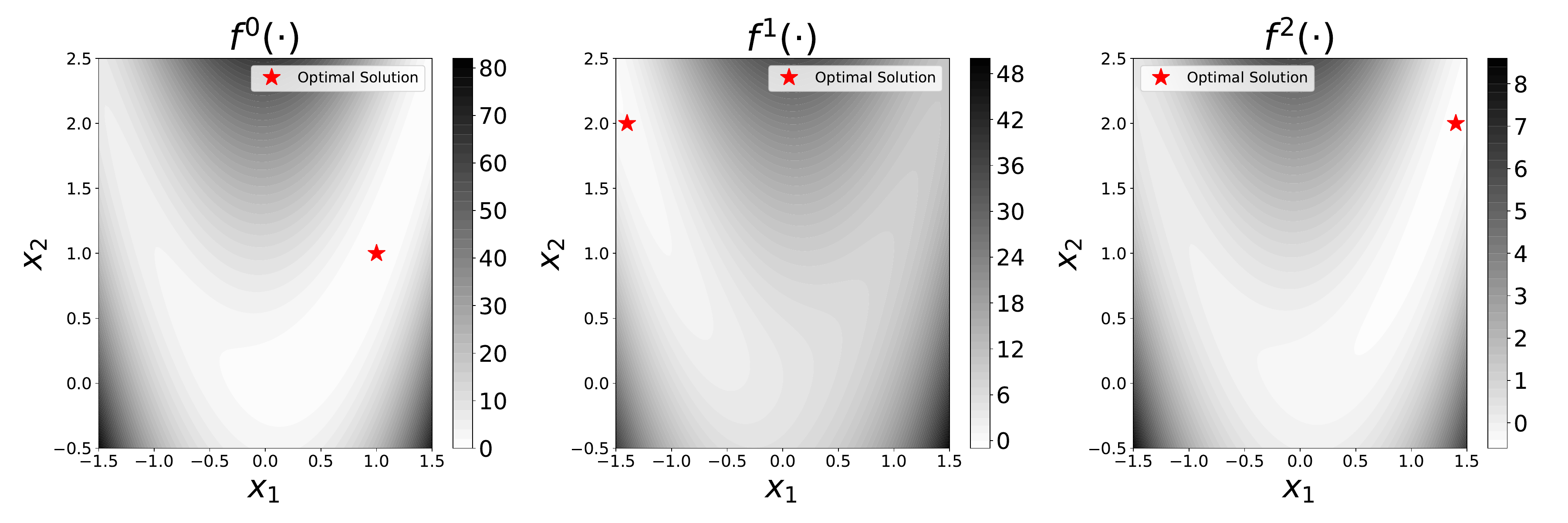}
\caption{Contour map of the deterministic Rosenbrock MF function for $d = 2$. The optimal solution is $(1, 1)$ for $f^0$, while the approximate optima for the lower-fidelity functions are $(-1.4, 2)$ for $f^1$ and $(1.4, 2)$ for $f^2$ with a box constraint $-2\le x_i \le 2$ for all $i\in\{1,2\}$.}
\label{fig:rosen-loss}
\end{figure}

The deterministic MF Rosenbrock functions were introduced in \cite{mainini2022mfexample} as $f^0(\BFx) = \sum_{i=1}^d 10(x_{i+1}-x_i^2)^2 + (1-x_i)^2,$ $f^1(\BFx) = \sum_{i=1}^{d-1} 5(x_{i+1}-x_i^2)^2 + (-2-x_i)^2 - \sum_{i=1}^d 0.5x_i$, and $f^2(\BFx) = (f^0-4-\sum_{i=1}^d0.5x_i)(10+\sum_{i=1}^d 0.25x_1)^{-1}.$ See Figure \ref{fig:rosen-loss} for two-dimensional loss landscapes. To introduce stochasticity, noise terms $E^t_i$ (drawn from a normal distribution with zero mean and dimension-dependent variance for any $i\in\{1,2,\dots,d\}$ and $t\in\{0,1,2\}$) are added: $\sum_{i=1}^d E_i^0$ to $f^0$, $\sum_{i=1}^d (E_i^0 + E_i^1)/2$ to $f^1$, and $\sum_{i=1}^d (E_i^0 + E_i^2)/2$ to $f^2$. The cost of querying each fidelity oracle is given by $\bm{w} = (1,\ 0.3,\ 0.1)$.

\begin{figure} [htp]
\centering
\subfloat[$\{\BFX_k\}$ trajectory with ASTRO-DF]{%
\resizebox*{7.5cm}{!}{\includegraphics{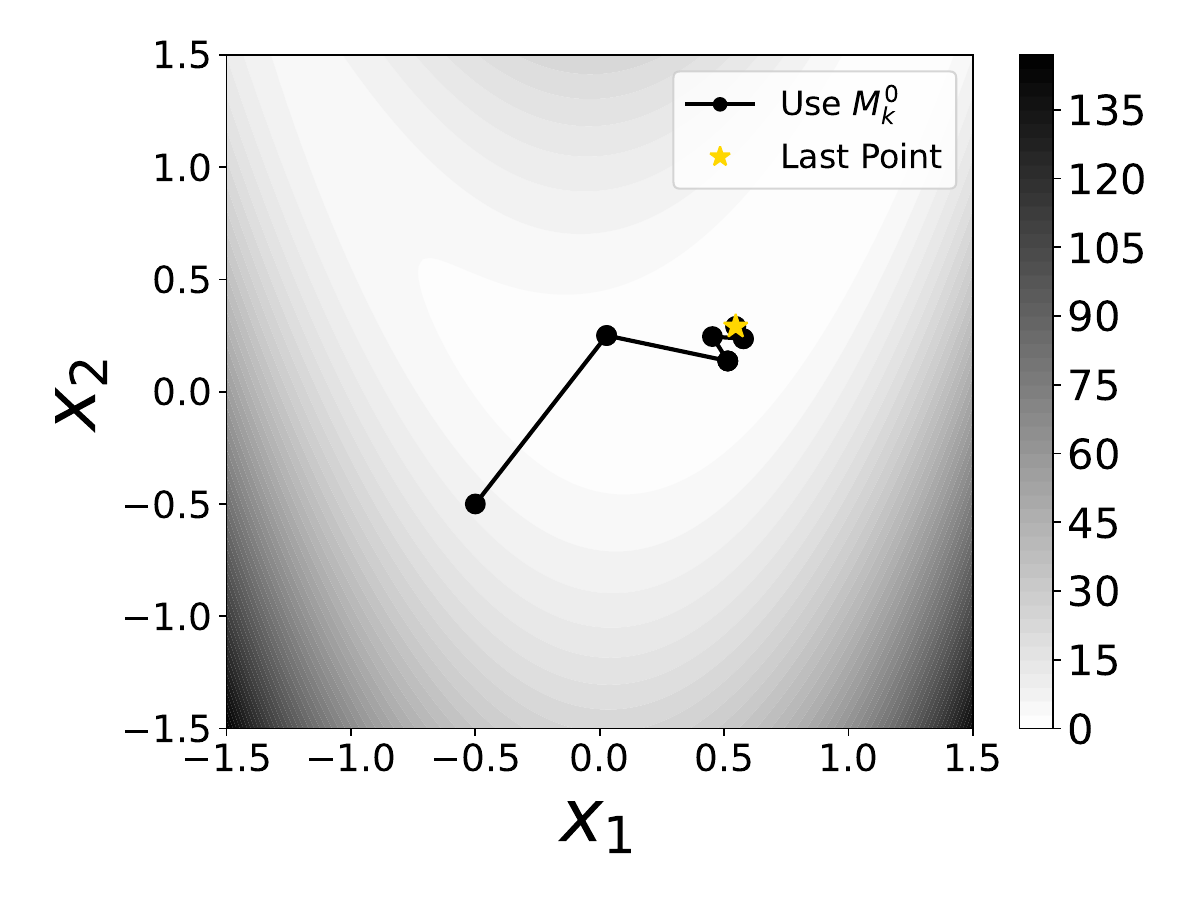}}\label{fig:astrodf-rosen}}
\subfloat[$\{\BFX_k\}$ trajectory with ASTRO-MFDF]{%
\resizebox*{7.5cm}{!}{\includegraphics{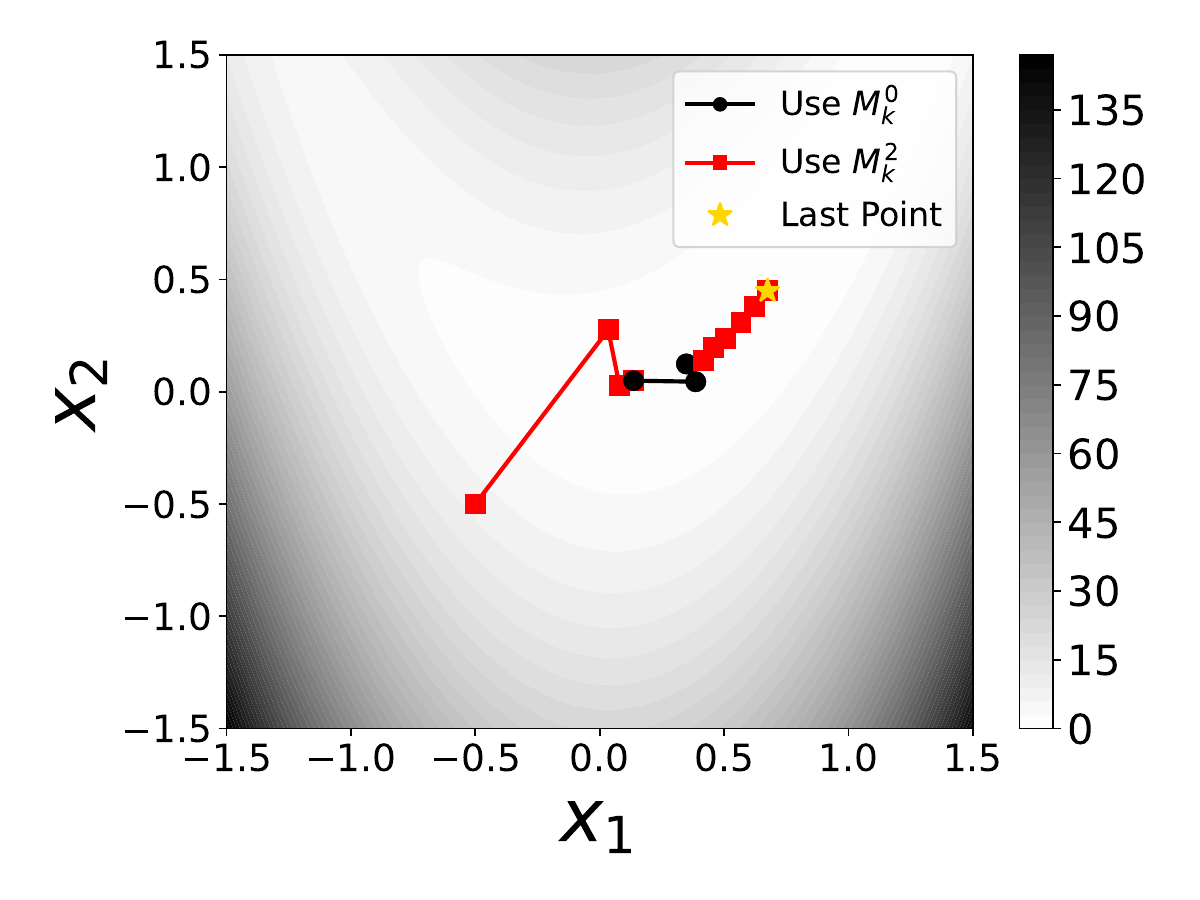}}\label{fig:astromfdf-rosen}}

\caption{One sample path of $\{\BFX_k\}$ with ASTRO-DF and ASTRO-MFDF on the 2-dimensional stochastic Rosenbrock function. Starting from the initial point (-0.5, -0.5) with a budget of 500 $0$-th fidelity oracle calls, the sequence $\{\BFX_k\}$ converges to (0.55, 0.29) in (a), with a corresponding $f^0$ value of approximately 0.308, and to (0.67, 0.45) in (b), with a corresponding $f^0$ value of approximately 0.109.} \label{fig:RSBR-trajectory}
\end{figure}

\begin{figure} [htp]
\centering
\subfloat[$d=10$]{%
\resizebox*{7.5cm}{!}{\includegraphics{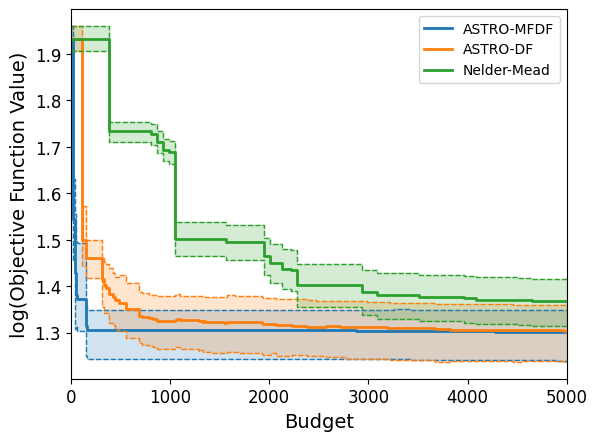}}\label{fig:rosen_10d}}
\subfloat[$d=20$]{%
\resizebox*{7.5cm}{!}{\includegraphics{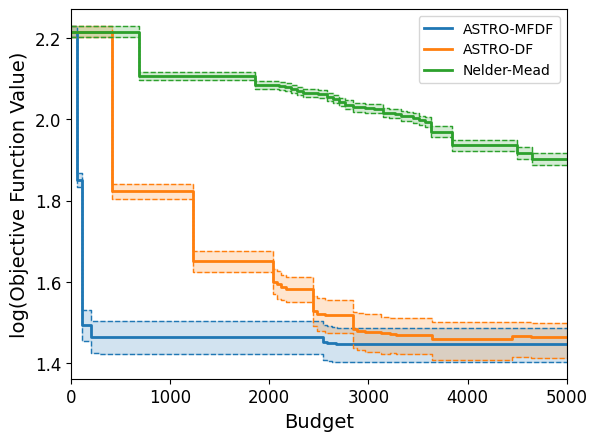}}\label{fig:rosen_20d}}
\caption{Finite-time performance on the stochastic Rosenbrock function, with a 95\% confidence interval with initial design point  $(-0.5) \times d$.  The x-axis represents a budget of 5000 $0$-th fidelity oracle calls, and the y-axis shows the objective function value on a logarithmic scale.} \label{fig:RSBR-highD}
\end{figure}

To illustrate the mechanism of ASTRO-MFDF, we compare its optimization trajectory with that of ASTRO-DF on the 2-dimensional stochastic Rosenbrock function (Figure \ref{fig:RSBR-trajectory}). As shown in Figure \ref{fig:astromfdf-rosen}, the optimization with ASTRO-MFDF begins using the lowest-fidelity function $f^2$, indicated in red. After several iterations, the local models constructed with lower-fidelity functions ($M^1_k$ and $M^2_k$) fail to yield improved iterates, prompting the algorithm to rely more heavily on the highest-fidelity function, represented in black. Later in the process, $M^2_k$ starts to produce better solutions again (red) enabling the algorithm to reduce computational cost per iteration. As a result, ASTRO-MFDF completes 24 iterations within the 500 $0$-th fidelity budget, compared to just 11 iterations by ASTRO-DF. It is worth noting that $M_k^1$ has not been used, as the correlation between $f^0$ and $f^1$ is expected to be low—i.e., $\alpha_k^1$ remains small—along the sample path of $\{\BFX_k\}$ generated by ASTRO-MFDF (see the contour map of $f^0$ and $f^1$ in Figure \ref{fig:rosen-loss}). This illustrates how ASTRO-MFDF can dynamically select the appropriate fidelity level at each iterate to maximize the computational efficiency of the MF approach. We also tested higher-dimensional problems, as shown in Figure \ref{fig:RSBR-highD}. As the dimensionality increases, ASTRO-MFDF demonstrates faster convergence, since constructing local models using lower-fidelity functions becomes significantly more cost-effective than using high-fidelity functions in high-dimensional settings.

\subsection{Experiments on continuous inventory problem}
In this section, we consider the $(s,S)$ inventory problem with continuous decision variables to create a more realistic experimental setting. The objective is to determine the optimal values of $s$ and $S$ that minimize the expected total cost, which includes holding costs, ordering costs, and backorder costs. Uncertainty in the system arises from two sources: (1) demand in each period follows an exponential distribution with mean $\theta$, and (2) lead times are drawn from a Poisson distribution with mean $\ell$ periods. For further details, see \cite{simoptgithub}. In the MF simulation, the 0-fidelity model runs for 100 days, the 1-fidelity model for 50 days, and the 2-fidelity model for 30 days, which implies the cost ratio vector $\bm{w}=(1,0.5,0.3)$.


\begin{figure} [htp]
\centering
\includegraphics[width=1\columnwidth]{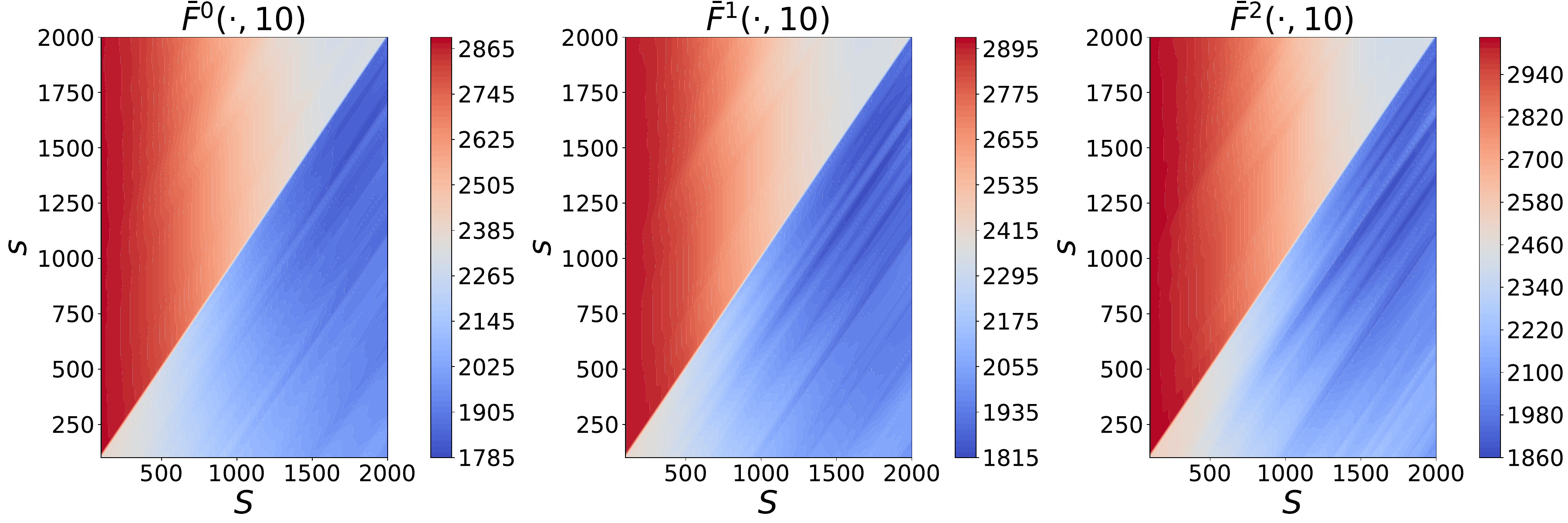}
\caption{Contour map of the estimated expected total cost for the $(s,S)$ inventory problem with $\theta = 400$ and $\ell=3$. Since the true objective function is unknown, the plots are based on estimates obtained from 10 samples.}
\label{fig:sscont-loss}
\end{figure}

\begin{figure} [htp]
\centering
\subfloat[$\{\BFX_k\}$ trajectory with ASTRO-DF]{%
\resizebox*{7.5cm}{!}{\includegraphics{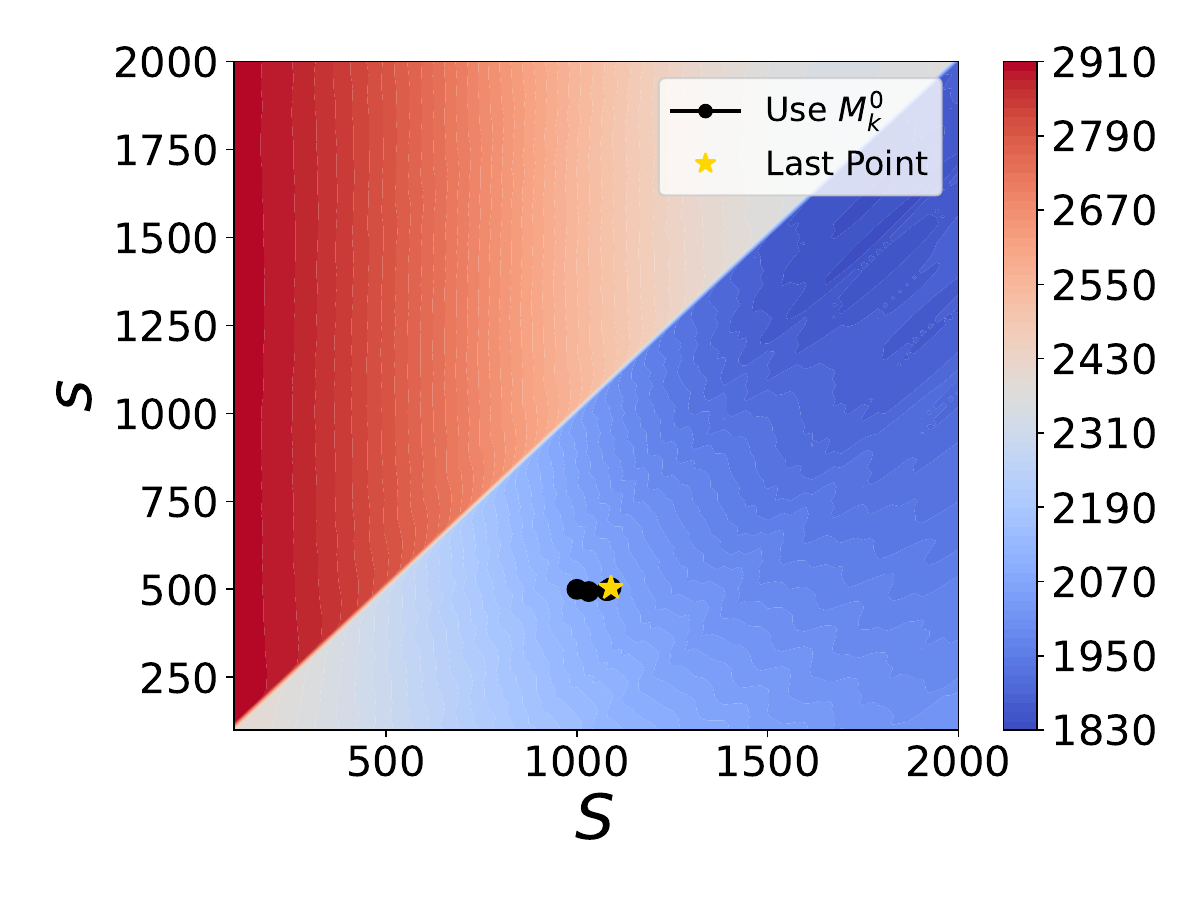}}\label{fig:solv-10}}
\subfloat[$\{\BFX_k\}$ trajectory with ASTRO-MFDF]{%
\resizebox*{7.5cm}{!}{\includegraphics{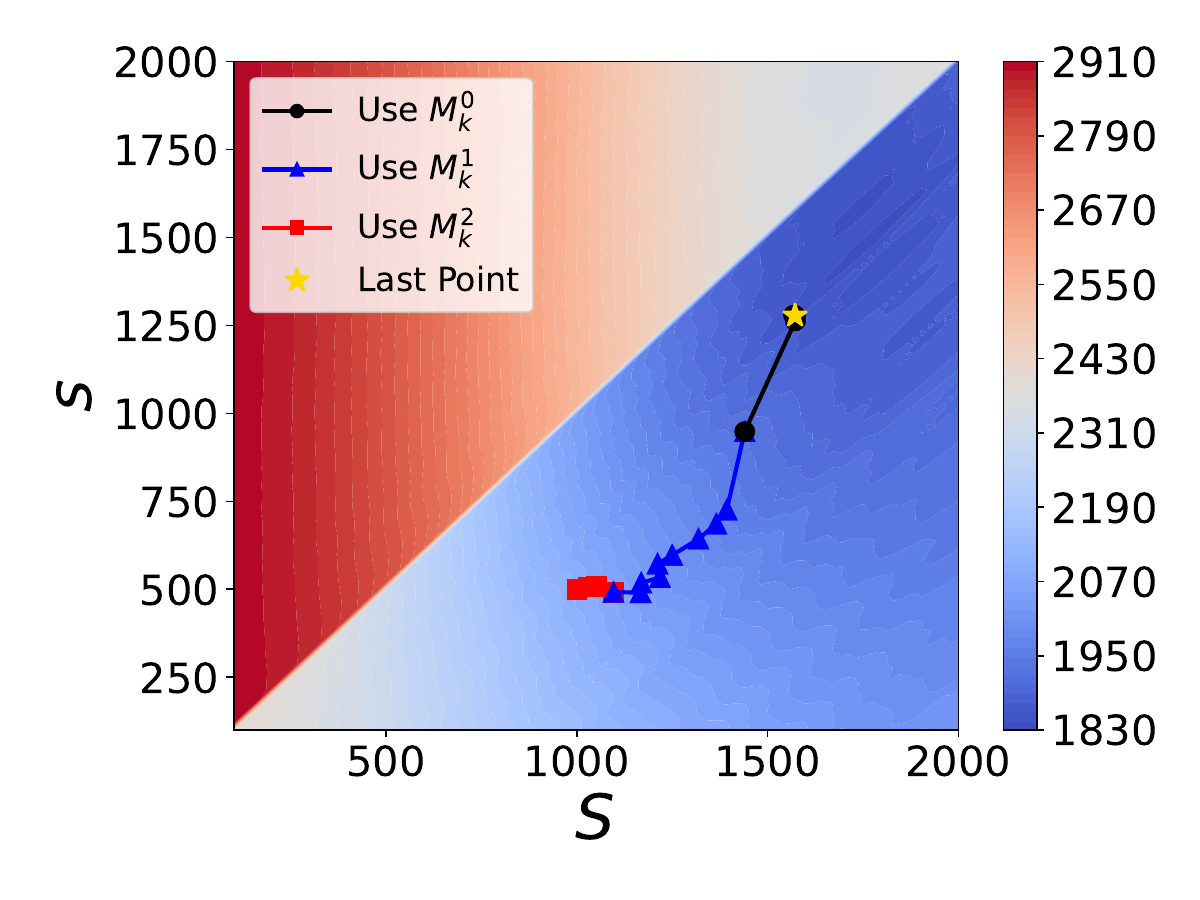}}\label{fig:solv-1}}
\caption{One sample path of $\{\BFX_k\}$ with ASTRO-DF and ASTRO-MFDF on the continuous $(s,S)$ inventory problem with $\theta = 400$ and $\ell = 3$. Starting from the initial point (500, 1000) with a budget of 1000 oracle calls, the sequence $\{\BFX_k\}$ converges to (504.13, 1089.27) in (a), and to (1277.85, 1571.87) in (b). The contour map shows the estimated expected total cost with 100 samples.} \label{fig:sscont-trajectory}
\end{figure}

\begin{figure} [htp]
\centering
\subfloat[$\BFx_0 = (200,500)$]{%
\resizebox*{7.5cm}{!}{\includegraphics{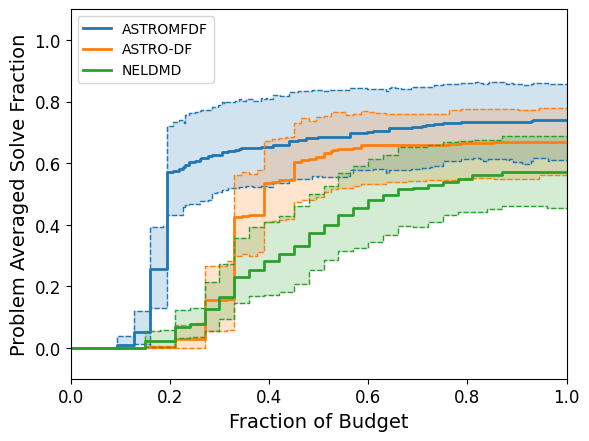}}\label{fig:solv-10}}
\subfloat[$\BFx_0 = (500,1000)$]{%
\resizebox*{7.5cm}{!}{\includegraphics{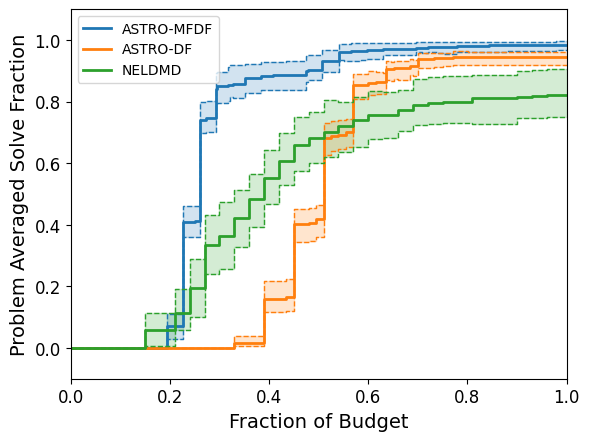}}\label{fig:solv-1}}
\caption{The fraction of problems each method successfully solved under 10\% optimality gap threshold.} \label{fig:solvability}
\end{figure}

As shown in Figure \ref{fig:sscont-loss}, the estimated objective function with a small sample size appears highly non-convex and non-smooth due to stochastic noise. As a result, the iterates are more likely to get stuck in regions that may not even correspond to local minima of the true objective function. When this happens, it typically indicates that the step size or search space has already become too small, making it difficult to escape without a significant amount of additional computational effort. As it can be seen in Figure \ref{fig:sscont-trajectory}, ASTRO-MFDF naturally addresses this issue by utilizing lower-fidelity models to escape suboptimal regions and preserving a large trust region for the high-fidelity function, thereby minimizing unnecessary computational costs. We also explored a range of settings for $\theta$ and $\ell$, with $\theta \in \{25, 50, 75, 100\}$ and $\ell \in \{1, 2, 3, 4, 5\}$. As shown in Figure \ref{fig:solvability}, ASTRO-MFDF achieves faster convergence on most problems and, in several cases, identifies better solutions, compared to ASTRO-DF and Nelder-Mead, in steady-state simulation optimization.

\section{Conclusion}
In MF methods, information from lower-fidelity models can sometimes hinder, rather than help, the optimization process. Therefore, it is crucial to sequentially determine whether to incorporate lower-fidelity information based on the current state of the optimization. Without this selective approach, using multiple information sources may actually slow convergence. To address this challenge, we introduce a novel stochastic trust-region method, ASTRO-MFDF, designed for MFSO. As a key feature of ASTRO-MFDF, a new adaptive sampling strategy, MFAS, is proposed utilizing MFMC to reduce the variance of function estimates. MFAS dynamically determines the sample sizes for the MF simulations and chooses between MFMC and MC based on sequential estimates of variance and covariance. MFAS helps minimize the waste of computational resources in several ways, such as reducing variance in function estimates and updating correlation vectors when reusing simulation outputs from past runs. Another important feature of ASTRO-MFDF is the incorporation of a correlation vector, which is updated dynamically using information from the optimization history. In particular, if the lower-fidelity models have not contributed to better solutions in past iterations, the method shifts focus to primarily use high-fidelity simulations for the optimization. Through numerical experiments with stochastic Rosenbrock functions and continuous $(s,S)$ inventory problems, we demonstrated that ASTRO-MFDF can achieve faster convergence. This is largely due to its ability to preserve a large trust region for high-fidelity simulations while efficiently leveraging low-fidelity ones, enabling more iterations under a limited budget. {For future research, we plan to address high-dimensional and computationally challenging traffic signal control problems in real-world networks by leveraging ASTRO-MFDF with subspace methods and high-performance computing.}

\section*{ACKNOWLEDGMENTS}
This work was authored in part by the National Renewable Energy Laboratory, operated by Alliance for Sustainable Energy, LLC, for the U.S. Department of Energy (DOE) under Contract No. DE-AC36-08GO28308. Funding was provided by the Office of Science, Office of Advanced Scientific Computing Research, Scientific Discovery through Advanced Computing (SciDAC) program through the FASTMath Institute. The views expressed in the article do not necessarily represent the views of the DOE or the U.S. Government. The U.S. Government retains and the publisher, by accepting the article for publication, acknowledges that the U.S. Government retains a nonexclusive, paid-up, irrevocable, worldwide license to publish or reproduce the published form of this work, or allow others to do so, for U.S. Government purposes.

\bibliography{main-paper}   

\begin{thebibliography}{10}

\bibitem{ha2025two}
Y.~Ha, S.~Shashaani, and M.~Menickelly, ``Two-stage estimation and variance modeling for latency-constrained variational quantum algorithms,'' {\em INFORMS Journal on Computing}, vol.~37, no.~1, pp.~125--145, 2025.

\bibitem{sakki2022renewable}
G.~K. Sakki, I.~Tsoukalas, P.~Kossieris, C.~Makropoulos, and A.~Efstratiadis, ``Stochastic simulation-optimization framework for the design and assessment of renewable energy systems under uncertainty,'' {\em Renewable and Sustainable Energy Reviews}, vol.~168, p.~112886, 2022.

\bibitem{zhang2022improved}
Z.~Zhang, Z.~Guan, Y.~Gong, D.~Luo, and L.~Yue, ``Improved multi-fidelity simulation-based optimisation: application in a digital twin shop floor,'' {\em International Journal of Production Research}, vol.~60, no.~3, pp.~1016--1035, 2022.

\bibitem{chen2017stochastic}
X.~Chen, S.~Hemmati, and F.~Yang, ``Stochastic co-kriging for steady-state simulation metamodeling,'' in {\em 2017 Winter Simulation Conference (WSC)} (W.~K.~V. Chan, A.~D'Ambrogio, G.~Zacharewicz, N.~Mustafee, G.~Wainer, and E.~Page, eds.), pp.~1750--1761, IEEE, 2017.

\bibitem{do2023multi}
B.~Do and R.~Zhang, ``Multi-fidelity bayesian optimization: A review,'' {\em arXiv preprint arXiv:2311.13050}, 2024.

\bibitem{foumani2023effects}
Z.~Z. Foumani, A.~Yousefpour, M.~Shishehbor, and R.~Bostanabad, ``On the effects of heterogeneous errors on multi-fidelity bayesian optimization,'' {\em arXiv preprint arXiv:2309.02771}, 2023.

\bibitem{Sara2018ASTRO}
S.~Shashaani, F.~S~Hashemi, and R.~Pasupathy, ``{ASTRO-DF}: A class of adaptive sampling trust-region algorithms for derivative-free stochastic optimization,'' {\em SIAM Journal on Optimization}, vol.~28, no.~4, pp.~3145--3176, 2018.

\bibitem{katya:DFObook}
A.~R. Conn, K.~Scheinberg, and L.~N. Vicente, {\em Introduction to derivative-free optimization}.
\newblock Society for Industrial and Applied Mathematics, 1st~ed., 2009.

\bibitem{ha2024adaptive}
Y.~Ha and J.~Mueller, ``Adaptive sampling-based bi-fidelity stochastic trust region method for derivative-free stochastic optimization,'' {\em arXiv preprint arXiv:2408.04625}, 2024.

\bibitem{karen2016mfmc}
B.~Peherstorfer, K.~Willcox, and M.~Gunzburger, ``Optimal model management for multifidelity monte carlo estimation,'' {\em SIAM Journal on Scientific Computing}, vol.~38, no.~5, pp.~A3163--A3194, 2016.

\bibitem{ha2023}
Y.~Ha, S.~Shashaani, and R.~Pasupathy, ``Complexity of zeroth-and first-order stochastic trust-region algorithms,'' {\em arXiv preprint arXiv:2405.20116}, 2024.

\bibitem{STRONG}
K.-H. Chang, L.~J. Hong, and H.~Wan, ``Stochastic trust-region response-surface method (strong)—a new response-surface framework for simulation optimization,'' {\em INFORMS Journal on Computing}, vol.~25, no.~2, pp.~230--243, 2013.

\bibitem{chen2018storm}
R.~Chen, M.~Menickelly, and K.~Scheinberg, ``Stochastic optimization using a trust-region method and random models,'' {\em Mathematical Programming}, vol.~169, no.~2, pp.~447--487, 2018.

\bibitem{Nocedal:NumericalOptbook}
J.~Nocedal and S.~J. Wright, {\em Numerical Optimization}.
\newblock Springer New York, 2nd~ed., 2006.

\bibitem{eckman2023simopt2}
D.~J. Eckman, S.~G. Henderson, and S.~Shashaani, ``Simopt: A testbed for simulation-optimization experiments,'' {\em INFORMS Journal on Computing}, vol.~35, no.~2, pp.~495--508, 2023.

\bibitem{mainini2022mfexample}
L.~Mainini, A.~Serani, M.~P. Rumpfkeil, E.~Minisci, D.~Quagliarella, H.~Pehlivan, S.~Yildiz, S.~Ficini, R.~Pellegrini, F.~Di~Fiore, {\em et~al.}, ``Analytical benchmark problems for multifidelity optimization methods,'' {\em arXiv preprint arXiv:2204.07867}, 2022.

\bibitem{simoptgithub}
D.~J. Eckman, S.~G. Henderson, S.~Shashaani, and R.~Pasupathy, ``{SimOpt}.'' \url{https://github.com/simopt-admin/simopt}, 2023.

\end{thebibliography}

\end{document}